\documentstyle[11pt]{article}
\begin{document}
\setlength{\textheight}{574pt}
\setlength{\textwidth}{432pt}
\setlength{\oddsidemargin}{18pt}
\setlength{\topmargin}{14pt}
\setlength{\evensidemargin}{18pt}
\newtheorem{theorem}{Theorem}[section]
\newtheorem{lemma}{Lemma}[section]
\newtheorem{sublemma}{Sublemma}[section]
\newtheorem{corollary}{Corollary}[section]
\newtheorem{remark}{Remark}[section]
\newtheorem{definition}{Definition}[section]
\newtheorem{problem}{Problem}
\newtheorem{proposition}{Proposition}[section]
\title{{\bf PLURICANONICAL SYSTEMS OF PROJECTIVE 3-FOLD OF GENERAL TYPE}}
\date{March, 2002}
\author{Hajime TSUJI}
\maketitle
\begin{abstract}
We prove that for every positive integer $m\geq 18(2^{9}\cdot 3^{7})!$ 
and every smooth projective $3$-fold of 
general type $X$ defined over complex numbers, 
$\mid mK_{X}\mid$ gives a birational rational map from $X$
into a projective space.
\end{abstract}
\tableofcontents
\section{Introduction}

Let $X$ be a smooth projective variety and let $K_{X}$ be the canonical 
bundle of $X$. 
$X$ is said to be a general type., if there exists a positive 
integer $m$ such that the pluricanonical system 
$\mid mK_{X}\mid$ gives a birational (rational) embedding of $X$. 
The following problem is fundamental to study projective 
vareity of general type. \vspace{10mm}\\
{\bf Problem}
Let $X$ be a smooth projective variety of general type.
Find a positive integer $m_{0}$ such that for every $m\geq m_{0}$,
$\mid mK_{X}\mid$ gives a birational rational map from $X$
into a projective space. \vspace{10mm} \\
If $\dim X = 1$, it is well known that $\mid 3K_{X}\mid$ gives a 
projective embedding.
In the case of smooth projective surfaces of general type, 
E. Bombieri showed that $\mid 5K_{X}\mid$ gives a birational
rational map from $X$ into a projective space (\cite{b3}).
In the case of $\dim X \geq 3$, recently the author 
gave a proof that there exists a positive integer $\nu_{n}$
depending only on $n$ 
such that for every smooth projective $n$-fold $X$ of general 
type, $\mid mK_{X}\mid$ gives a birational embedding of $X$
into a projective space for every $m\geq \nu_{n}$
(\cite[Theorem 1.1]{tu5}).

But even in  the case of 3-fold, it has not yet been known
the explicit upper bound for $\nu_{3}$. 
In this article I would like to give an explicit upper bound 
for $\nu_{3}$. 

\begin{theorem} 
For every posive integer $m \geq 18(2^{9}\cdot 3^{7})!$ and smooth projective $3$-fold $X$
of general type defined over complex numbers, $\mid mK_{X}\mid$ gives a birational rational map
from $X$ into a projective space.
\end{theorem}
The bound here is  astronomous and far from being optimal. 
Yes, this bound is definitely nonsense. 
But this is the first explicit bound for $\nu_{3}$.
I am just interested in giving an explicit bound. 

It is not difficult to improve the bound for $\nu_{3}$ 
given in Theorem 1.1.  
In fact there are a lot of obvious losses in the 
estimates given in this paper. 
But I decided not to include such an improvement, 
because I think it is completely nonsense. 
 
One reason is that  as is shown in \cite{r}, the actual bound 
for $\nu_{3}$ may be unexpectedly large. 
Actually in \cite[p.359, see $X_{46}$]{r} M. Reid constructed a minimal 3-fold $X$ of general type 
with $K_{X}^{3} = 1/420$ (this $X$ is a complete intersection 
in a weighted projective space). 
This already implies that $\nu_{3} \geq 9$.

Another reason is that it is very plausible that we may encounter 
a very strange 3-fold of general type near future.
In my opinion, even if there is  some hope 
to find the best bound for $\nu_{3}$,
this kind of result may be uninteresting, 
unless it gives a perspective to the best bounds 
for $\nu_{n}$ for all $n$, because I think that  
to consider 3-folds only is meaningless. 

We note that there are several previous works \cite{be,ma,l,l2,c1,c2} 
by X. Benveniste, K. Matuki, T. Luo and M. Chen. 
Except T. Luo, they put conditions on minimal models.
But since the construction of minimal models for 3-folds
is very far from being explicit, I think that   
if we put condition on minimal models,
the meaning of the results will be very much destroyed.

In contrast to the theory of curves and surfaces, 
the higher dimensional algebraic geometry 
is  rather philosophical than 
practical or concrete. 
And we should give up very refined and 
concrete results and  change our focus to more abstract 
theorems,   
otherwise we will be involved in tremendously complicated 
situation and will be very hard to continue. 

This paper consists of two parts. 

The first part is a review of the arguments in \cite{tu5}. 
The argument here works in all dimensions, if we assume 
the minimal model programs. 

The second part is special for 3-folds of general type. 
Here we use a detailed description of 3-dimensional 
terminal singularities (\cite{m}).

So far it seems to be hopeless to get a perspective 
about explicit bounds for $\nu_{n}$ for $n\geq 4$. 

In this paper all the varieties are defined over {\bf C}.

\section{How to bound the degree} 

In this section, we shall prove the following weaker theorem 
which is the special case of Theorem 1.1 in \cite{tu5}. 

\begin{theorem}(\cite[Theorem 1.1]{tu5})
There exists a positive integer $\nu_{3}$ such that 
for every projective 3-fold of general type $X$ and 
for every $m\geq \nu_{3}$, 
$\mid mK_{X}\mid$ gives a birational embedding of $X$ 
into a projective space. 
\end{theorem}

Since we need not only the result itself but also 
the proof, here we include the full proof of 
Theorem 2.1. 
The proof presented here is a part of the proof in \cite{tu5},
but it is much easier, since there exists a minimal model 
for a projective 3-fold of general type. 

Let $X$ be a minimal projective $3$-fold of general type, i.e.,
$X$ has only {\bf Q}-factorial terminal singularities and the 
canonical divisor $K_{X}$ is nef. 
We set 
\[
X^{\circ} = \{ x\in X\mid x\not{\in} \mbox{Bs}\mid mK_{X}\mid \mbox{and  
$\Phi_{\mid mK_{X}\mid}$ is a biholomorphism} 
\]
\[
\hspace{50mm} \mbox{on a neighbourhood of $x$}\} .
\]
Then $X^{\circ}$ is a nonempty Zariski open subset of $X$. 
 
\subsection{Construction of a stratification}
We set 
\[
\mu_{0}: = K_{X}^{3}.
\]
\begin{lemma} Let $x,x^{\prime}$ be distinct points on $X^{\circ}$.  
We set 
\[
{\cal M}_{x,x^{\prime}} = {\cal M}_{x}\otimes{\cal M}_{x^{\prime}},
\]
where ${\cal M}_{x},{\cal M}_{x^{\prime}}$ denote the
maximal ideal sheaf of the points $x,x^{\prime}$ respectively.
Let $\varepsilon$ be a sufficiently small positive number.
Then 
\[
H^{0}(X,{\cal O}_{X}(mK_{X})\otimes{\cal M}_{x,x^{\prime}}^{\lceil\sqrt[3]{\mu_{0}}
(1-\varepsilon )\frac{m}{\sqrt[3]{2}}\rceil})\neq 0
\]
for every sufficiently large $m$.
\end{lemma}
{\bf Proof of Lemma 2.1}.   
Let us consider the exact sequence:
\[
0\rightarrow H^{0}(X,{\cal O}_{X}(mK_{X})\otimes
{\cal M}_{x,x^{\prime}}^{\lceil\sqrt[3]{\mu_{0}}(1-\varepsilon )\frac{m}{\sqrt[n]{2}}\rceil})
\rightarrow H^{0}(X,{\cal O}_{X}(mK_{X}))\rightarrow
\]
\[
  H^{0}(X,{\cal O}_{X}
(mK_{X})/{\cal M}_{x,x^{\prime}}^{\lceil\sqrt[3]{\mu_{0}}(1-\varepsilon )\frac{m}{\sqrt[3]{2}}\rceil}).
\]
We note that 
\[
3!\cdot\overline{\lim}_{m\rightarrow\infty}m^{-3}\dim H^{0}(X,{\cal O}_{X}(mK_{X})) = \mu_{0}
\]
holds, since $K_{X}$ is nef and big. 
Since
\[
3!\cdot\overline{\lim}_{m\rightarrow\infty}m^{-3}\dim H^{0}(X,{\cal O}_{X}(mK_{X})
/{\cal M}_{x,x^{\prime}}^{\lceil\sqrt[3]{\mu_{0}}(1-\varepsilon )\frac{m}{\sqrt[3]{2}}\rceil})
=
\mu_{0}(1-\varepsilon )^{3} < \mu_{0}
\]
hold, we see that Lemma 2.1 holds.  {\bf Q.E.D.}

\vspace{5mm}

Let us take a sufficiently large positive integer $m_{0}$ and let $\sigma$
be a general (nonzero) element $\sigma_{0}$ of  
$H^{0}(X,{\cal O}_{X}(m_{0}K_{X})\otimes
{\cal M}_{x,x^{\prime}}^{\lceil\sqrt[3]{\mu_{0}}(1-\varepsilon )\frac{m_{0}}{\sqrt[3]{2}}\rceil})$. 
We define an effective {\bf Q}-divisor $D_{0}$ by 
\[
D_{0} = \frac{1}{m_{0}}(\sigma_{0}).
\]
We define a positive number $\alpha_{0}$ by
\[
\alpha_{0} := \inf\{\alpha > 0\mid 
\mbox{$(X,\alpha D_{0})$ is not KLT at 
both $x$ and $x^{\prime}$}\},
\]
where KLT is short for of Kawamata logterminal (cf. \cite[p.56, Definition 2.34]{k-m}). 
Since $(\sum_{i=1}^{3}\mid z_{i}\mid^{2})^{-3}$ is not locally integrable 
around $O\in \mbox{{\bf C}}^{3}$, by the construction of $D_{0}$, we see
that 
\[
\alpha_{0}\leq \frac{3\sqrt[3]{2}}{\sqrt[3]{\mu_{0}}(1-\varepsilon )}
\]
holds.
Then one of the following two cases occurs. \vspace{5mm} \\
{\bf Case} 1.1:  For every small positive number $\delta$, 
$(X,(\alpha_{0}-\delta )D_{0})$ is KLT 
at both $x$ and $x^{\prime}$. \\
{\bf Case} 1.2: For every small positive number $\delta$, 
$(X,(\alpha_{0}-\delta )D_{0})$ is KLT  
at one of $x$ or $x^{\prime}$ say $x$. \vspace{10mm} \\

We first consider Case 1.1.
Let $X_{1}$ be the minimal center of log canonical singularities 
at $x$ (cf. \cite{ka2}). 
We consider the following two cases. \vspace{5mm} \\
{\bf Case} 2.1: $X_{1}$ passes through both $x$ and $x^{\prime}$, \\
{\bf Case} 2.2: Otherwise \vspace{10mm} \\

For the first we consider Case 2.1.
In this case $X_{1}$ is not isolated at $x$.  Let $n_{1}$ denote
the dimension of $X_{1}$.  Let us define the volume $\mu_{1}$ of $X_{1}$
with respect to $K_{X}$ by
\[
\mu_{1} := K_{X}^{n_{1}}\cdot X_{1}.
\]
Since $x\in X^{\circ}$, we see that $\mu_{1} > 0$ holds.
The proof of the following lemma is identical that of Lemma 2.1.
\begin{lemma} Let $\varepsilon$ be a sufficiently small positive number and let $x_{1},x_{2}$ be distinct regular points on $X_{1}$. 
Then for a sufficiently large $m >1$,
\[
H^{0}(X_{1},{\cal O}_{X_{1}}(mK_{X})\otimes
{\cal M}_{x_{1},x_{2}}^{\lceil\sqrt[n_{1}]{\mu_{1}}(1-\varepsilon )\frac{m}{\sqrt[n_{1}]{2}}\rceil})\neq 0
\]
holds.
\end{lemma}
Let $x_{1},x_{2}$ be two distinct regular points on $X_{1}$. 
Let $m_{1}$ be a sufficiently large positive integer and 
Let 
\[
\sigma_{1}^{\prime}
\in 
H^{0}(X_{1},{\cal O}_{X_{1}}(mK_{X})\otimes
{\cal M}_{x_{1},x_{2}}^{\lceil\sqrt[n_{1}]{\mu_{1}}(1-\varepsilon )\frac{m}{\sqrt[n_{1}]{2}}\rceil})
\]
be a nonzero element. 

By Kodaira's lemma there is an effective {\bf Q}-divisor $E$ such
that $K_{X}- E$ is ample.
Let $\ell_{1}$ be a sufficiently large positive integer which will be specified later  such that
\[
L_{1} := \ell_{1}(K_{X}- E)
\]
is Cartier. 

\begin{lemma}
If we take $\ell_{1}$ sufficiently large, then 
\[
\phi_{m} : H^{0}(X,{\cal O}_{X}(mK_{X}+L_{1}))\rightarrow 
H^{0}(X_{1},{\cal O}_{X_{1}}(mK_{X}+L_{1} ))
\]
is surjective for every  $m\geq 0$.
\end{lemma}
{\bf Proof}.
$K_{X}$ is nef {\bf Q}-Cartier divisor by the assumption.
Let $r$ be the index of $X$, i.e. $r$ is the minimal positive integer such that $rK_{X}$ is Cartier.  Then since $K_{X}$ is semiample 
(\cite{ka0}), by the Kodaira-Nakano vanishing theorem, 
for every locally free sheaf 
${\cal E}$, 
there exists a positive integer $k_{0}$ 
depending on ${\cal E}$
such that for every $\ell \geq k_{0}$ 
\[
H^{q}(X,{\cal O}_{X}((1+mr)K_{X}+L_{1})\otimes {\cal E}) = 0
\]
holds for every $q\geq 1$ and $m\geq 0$.
Let us consider the exact sequences 
\[
0 \rightarrow {\cal K}_{j}\rightarrow {\cal E}_{j}
\rightarrow {\cal O}_{X}(jK_{X})\otimes{\cal I}_{X_{1}}\rightarrow 0 
\]
for some locally free sheaf ${\cal E}_{j}$ for every 
$0 \leq j \leq r-1$,
where ${\cal I}_{X_{1}}$ denotes the ideal sheaf associated with
$X_{1}$. 
Then noting the above fact, we can prove that 
if  we take $\ell_{1}$ sufficiently large,
\[
H^{q}(X,{\cal O}_{X}(mK_{X}+L_{1})\otimes {\cal I}_{X_{1}})
\]
holds for every $q\geq 1$ and $m\geq 0$
by exactly the same manner as the standard proof of 
Serre's vanishing theorem (cf. \cite[p.228, Theorem 5.2]{ha}).
This implies the desired surjection. 
\vspace{5mm}{\bf Q.E.D.} \\ 

Let $\tau$ be a general section in 
$H^{0}(X,{\cal O}_{X}(L_{1}))$.
Then by Lemma 2.3 we see that   
\[
\sigma_{1}^{\prime}\otimes\tau\in
H^{0}(X_{1},{\cal O}_{X_{1}}(m_{1}K_{X}+L_{1})
{\cal M}_{x_{1},x_{2}}^{\lceil\sqrt[n_{1}]{\mu_{1}}(1-\varepsilon )\frac{m_{1}}
{\sqrt[n_{1}]{2}}\rceil})
\]
extends to a section
\[
\sigma_{1}\in H^{0}(X,{\cal O}_{X}((m_{1}+\ell_{1} )K_{X})).
\]
We may assume that  there exists a neighbourhood $U_{x,x^{\prime}}$ of $\{ x,x^{\prime}\}$ such that the divisor $(\sigma _{1})$  is smooth
on  $U_{x,x^{\prime}} - X_{1}$ by Bertini's theorem, if we take $\ell_{1}$
sufficiently large, since as Lemma 2.3 
\[
H^{0}(X,{\cal O}_{X}(mK_{X}+L_{1})\otimes{\cal I}(h^{m}))
\rightarrow
H^{0}(X,{\cal O}_{X}(mK_{X}+L_{1})\otimes{\cal I}_{X_{1}}\cdot{\cal M}_{y})
\]
is surjective for every $y\in X$ and
 $m\geq 0$.
We set 
\[
D_{1} = \frac{1}{m_{1}+\ell_{1}}(\sigma_{1}). 
\]
Suppose that $x,x^{\prime}$ are nonsingular points on $X_{1}$.
Then we set $x_{1} = x, x_{2} = x^{\prime}$.
Let $\varepsilon_{0}$ be a sufficiently small 
positive rational number and define $\alpha_{1}$ by
\[
\alpha_{1} := \inf\{\alpha > 0 \mid 
\mbox{ $(\alpha_{0}-\varepsilon_{0})D_{0} + \alpha D_{1}$ 
is not KLT at both $x$ and $x^{\prime}$} \}.
\]
Then we may define the proper subvariety $X_{2}$ of $X_{1}$ 
as a minimal center of logcanonical singularities as before. 

By Lemma 2.2 we may assume that we have taken $m_{1}$ so that  
\[
\frac{\ell_{1}}{m_{1}} \leq 
\varepsilon_{0}\frac{\sqrt[n_{1}]{\mu_{1}}}{n_{1}\sqrt[n_{1}]{2}}
\]
holds.
\begin{lemma}
\[
\alpha_{1}\leq \frac{n_{1}\sqrt[n_{1}]{2}}{\sqrt[n_{1}]{\mu_{1}}} 
+ O(\varepsilon _{0})
\]
holds.
\end{lemma}
To prove Lemma 2.4, we need the following elementary lemma.
\begin{lemma}(\cite[p.12, Lemma 6]{t})
Let $a,b$ be  positive numbers. Then
\[
\int_{0}^{1}\frac{r_{2}^{2n_{1}-1}}{(r_{1}^{2}+r_{2}^{2a})^{b}}
dr_{2}
=
r_{1}^{\frac{2n_{1}}{a}-2b}\int_{0}^{r_{1}^{-{2}{a}}}
\frac{r_{3}^{2n_{1}-1}}{(1 + r_{3}^{2a})^{b}}dr_{3}
\]
holds, where 
\[
r_{3} = r_{2}/r_{1}^{1/a}.
\]
\end{lemma}
{\bf Proof of Lemma 2.4.}
Let $(z_{1},\ldots ,z_{n})$ be a local coordinate on a 
neighbourhood $U$ of $x$ in $X$ such that 
\[
U \cap X_{1} = 
\{ q\in U\mid z_{n_{1}+1}(q) =\cdots = z_{n}(q)=0\} .
\] 
We set $r_{1} = (\sum_{i=n_{1}+1}^{n}\mid z_{1}\mid^{2})^{1/2}$ and 
$r_{2} = (\sum_{i=1}^{n_{1}}\mid z_{i}\mid^{2})^{1/2}$.
Then there exists a positive constant $C$ such that 
\[
\parallel\sigma_{1}\parallel^{2}\leq 
C(r_{1}^{2}+r_{2}^{2\lceil\sqrt[n_{1}]{\mu_{1}}(1-\varepsilon )\frac{m_{1}}
{\sqrt[n_{1}]{2}}\rceil})
\]
holds on a neighbourhood of $x$, 
where $\parallel\,\,\,\,\parallel$ denotes the norm with 
respect to $h_{X}^{m_{1}+\ell_{1}}$.
We note that there exists a positive integer $M$ such that 
\[
\parallel\sigma\parallel^{-2} = O(r_{1}^{-M})
\]
holds on a neighbourhood of the generic point of $U\cap X_{1}$,
where $\parallel\,\,\,\,\parallel$ denotes the norm with respect to 
$h_{X}^{m_{0}}$. 
Then by Lemma 2.5., we have the inequality 
\[
\alpha_{1} \leq (\frac{m_{1}+\ell_{1}}{m_{1}})\frac{n_{1}\sqrt[n_{1}]{2}}{\sqrt[n_{1}]{\mu_{1}}} 
+ O(\varepsilon _{0})
\] 
holds. 
By using the fact that 
\[
\frac{\ell_{1}}{m_{1}} \leq 
\varepsilon_{0}\frac{\sqrt[n_{1}]{\mu_{1}}}{n_{1}\sqrt[n_{1}]{2}}
\]
we obtain that 
\[
\alpha_{1}\leq \frac{n_{1}\sqrt[n_{1}]{2}}{\sqrt[n_{1}]{\mu_{1}}} 
+ O(\varepsilon _{0})
\]
holds.
{\bf Q.E.D.} \vspace{5mm} \\
If $x$ or $x^{\prime}$ is a singular point on $X_{1}$, we need the following lemma.
\begin{lemma}
Let $\varphi$ be a plurisubharmonic function on $\Delta^{n}\times{\Delta}$.
Let $\varphi_{t}(t\in\Delta )$ be the restriction of $\varphi$ on
$\Delta^{n}\times\{ t\}$.
Assume that $e^{-\varphi_{t}}$ does not belong to $L^{1}_{loc}(\Delta^{n},O)$
for every $t\in \Delta^{*}$.

Then $e^{-\varphi_{0}}$ is not locally integrable at $O\in\Delta^{n}$.
\end{lemma}
Lemma 2.6 is an immediate consequence of \cite[p.20, Theorem]{o-t}.
Using Lemma 2.6 and Lemma 2.5, we see that Lemma 2.4 holds
by letting $x_{1}\rightarrow x$ and $x_{2}\rightarrow x^{\prime}$.

\vspace{5mm}

Next we consider Case 1.2 and Case 2.2.  
In this case  for every  sufficiently small positive number $\delta$, 
$(X,(\alpha_{0}-\varepsilon_{0})D_{0}+(\alpha_{1}-\delta )D_{1})$
is KLT at $x$ and not KLT at $x^{\prime}$. 

In these cases, instead of Lemma 2.2, we use the following simpler lemma.

\begin{lemma} Let $\varepsilon$ be a sufficiently small positive number and let $x_{1}$ be a smooth point on $X_{1}$. 
Then for a sufficiently large $m >1$,
\[
H^{0}(X_{1},{\cal O}_{X_{1}}(mK_{X})\otimes
{\cal M}_{x_{1}}^{\lceil\sqrt[n_{1}]{\mu_{1}}(1-\varepsilon )m\rceil})\neq 0
\]
holds.
\end{lemma}

Then taking a general nonzero element $\sigma_{1}^{\prime}$ in
\[
H^{0}(X_{1},{\cal O}_{X_{1}}(m_{1}K_{X})\otimes{\cal I}(h^{m_{1}})\otimes
{\cal M}_{x_{1}}^{\lceil\sqrt[n_{1}]{\mu_{1}}(1-\varepsilon )m_{1}
\rceil}),
\]
for a sufficiently large $m_{1}$.
As in Case 1.1 and Case 2.1 we obtain a proper subvariety
$X_{2}$ in $X_{1}$ also in this case.

Inductively for distinct points $x,x^{\prime}\in X^{\circ}$, we construct a strictly decreasing
sequence of subvarieties
\[
X = X_{0}\supset X_{1}\supset \cdots \supset X_{r}\supset X_{r+1} =  x\,\,\mbox{or}\,\, x^{\prime}
\]
and invariants (depending on small positive  rational numbers $\varepsilon_{0},\ldots ,
\varepsilon_{r-1}$, large positive integers $m_{0},m_{1},\ldots ,m_{r}$, etc.) :
\[
\alpha_{0} ,\alpha_{1},\ldots ,\alpha_{r},
\]
\[
\mu_{0},\mu_{1},\ldots ,\mu_{r}
\]
and
\[
n >  n_{1}> \cdots > n_{r}.
\]
By Nadel's vanishing theorem (\cite[p.561]{n}) we have the following lemma.
\begin{lemma} 
Let $x,x^{\prime}$ be two distinct points on $X^{\circ}$. 
Then for every $m\geq \lceil\sum_{i=0}^{r}\alpha_{i}\rceil +1$,
$\Phi_{\mid mK_{X}\mid}$ separates $x$ and $x^{\prime}$.
\end{lemma}
{\bf Proof}. 
For $i= 0,1,\ldots, r$ let $h_{i}$ be the singular hermitian metric 
on $K_{X}$ defined by 
\[
h_{i}:= \frac{1}{\mid\sigma_{i}\mid^{\frac{2}{m_{i}+\ell_{i}}}},
\]
where we have set $\ell_{0}:= 0$. 
More precisely for any $C^{\infty}$-hermitian metric $h_{X}$ 
on $K_{X}$ we have defined $h_{i}$ as 
\[
h_{i}:= \frac{h_{X}}{h_{X}^{m_{i}+\ell_{i}}(\sigma_{i},\sigma_{i})^{\frac{1}{m_{i}+\ell_{i}}}}.
\]
Using Kodaira's lemma (\cite[Appendix]{k-o}), 
let $E$ be an effective {\bf Q}-divisor $E$ such that
$K_{X} - E$ is ample. 
Let $m$ be a positive integer such that $m\geq \lceil\sum_{i=0}^{r}\alpha_{i}\rceil +1$ holds. 
Let $h_{L}$ is a $C^{\infty}$-hermitian metric on the ample {\bf Q}-line bundle 
\[
L := (m-1-(\sum_{i=0}^{r-1}(\alpha_{i}-\varepsilon_{i}))
-\alpha_{r}-\delta_{L})K_{X} - \delta_{L}E
\]
 with strictly positive curvature, where $\delta_{L}$ be a sufficiently small positive number 
and we have considered $h_{L}$ as a singular hermitian metric 
on  \\ $(m-1- (\sum_{i=0}^{r-1}(\alpha_{i}-\varepsilon_{i}))-\alpha_{r})K_{X}$.
Let us define the singular hermitian metric $h_{x,x^{\prime}}$ of $(m-1)K_{X}$ defined by  
\[
h_{x,x^{\prime}} = (\prod_{i=0}^{r-1}h_{i}^{\alpha_{i}-\varepsilon_{i}})\cdot 
 h_{r}^{\alpha_{r}}\cdot h_{L}.
\]
Then we see that  ${\cal I}(h_{x,x^{\prime}})$ defines a subscheme of 
$X$ with isolated support around $x$ or $x^{\prime}$ by the definition of 
the invariants $\{\alpha_{i}\}$'s. 
By the construction the curvature current $\Theta_{h_{x,x^{\prime}}}$ is strictly positive on $X$. 
Then by Nadel's vanishing theorem (\cite[p.561]{n}) we see that 
\[
H^{1}(X,{\cal O}_{X}(mK_{X})\otimes {\cal I}(h_{x,x^{\prime}})) = 0
\]
holds. 
Hence 
\[
H^{0}(X,{\cal O}_{X}(mK_{X}))
\rightarrow 
H^{0}(X,{\cal O}_{X}\otimes {\cal O}_{X}/{\cal I}(h_{x,x^{\prime}}))
\]
is surjective. 
Since by the construction of $h_{x,x^{\prime}}$ (if we take
$\delta_{L}$ sufficiently small) 
$\mbox{Supp}({\cal O}_{X}/{\cal I}(h_{x,x^{\prime}}))$ 
contains both $x$ and $x^{\prime}$ and is 
isolated at least one of $x$ or $x^{\prime}$. 
Hence by the above surjection, there exists a section
$H^{0}(X,{\cal O}_{X}(mK_{X}))$ such that 
\[
\sigma (x) \neq 0,\sigma (x^{\prime}) = 0
\]
or 
\[
\sigma (x) = 0,\sigma (x^{\prime}) \neq 0
\]
holds. 
This implies that $\Phi_{\mid mK_{X}\mid}$ separates 
$x$ and $x^{\prime}$.   {\bf Q.E.D.} 

\subsection{Construction of the stratification as a family}

In this subsection we shall construct the above stratification as a family. 

We note that for a fixed pair $(x,x^{\prime}) \in X^{\circ}\times X^{\circ}-\Delta_{X}$, $\sum_{i=0}^{r}\alpha_{i}$ depends on the choice of $\{ X_{i}\}$'s, where 
$\Delta_{X}$ denotes the diagonal of $X\times X$. 
Moving $(x,x^{\prime})$ in $X^{\circ}\times X^{\circ} - \Delta_{X}$, 
we  shall consider the above operation simultaneously.
Let us explain the procedure. 
We set 
\[
B := X^{\circ}\times X^{\circ} - \Delta_{X}.
\] 
Let 
\[
p : X\times B\longrightarrow X
\]
be the first projection and let   
\[
q : X\times {B}
\longrightarrow B
\]
be the second projection. 
Let $Z$ be the subvariety of $X\times B$ defined by
\[
Z := \{ (x_{1},x_{2},x_{3}) : X\times B \mid 
x_{1} = x_{2} \,\,\mbox{or} \,\, x_{1} = x_{3} \} .
\]
In this case we consider 
\[
q_{*}{\cal O}_{X\times B}(m_{0}p^{*}K_{X})
\otimes {\cal I}_{Z}^{\lceil\sqrt[3]{\mu_{0}}(1-\varepsilon )\frac{m_{0}}{\sqrt[3]{2}}\rceil }
\]
instead of 
\[
H^{0}(X,{\cal O}_{X}(m_{0}K_{X})\otimes
{\cal M}_{x,x^{\prime}}^{\lceil\sqrt[3]{\mu_{0}}(1-\varepsilon )\frac{m_{0}}{\sqrt[3]{2}}\rceil}),
\]
where ${\cal I}_{Z}$ denotes the ideal sheaf of $Z$. 
Let $\tilde{\sigma}_{0}$ be a nonzero global meromorphic 
section  of 
\[
q_{*}{\cal O}_{X\times B}(m_{0}p^{*}K_{X})
\otimes {\cal I}_{Z}^{\lceil\sqrt[3]{\mu_{0}}(1-\varepsilon )\frac{m_{0}}{\sqrt[3]{2}}\rceil } 
\]
on $B$ for a sufficiently large positive integer $m_{0}$.
We set 
\[
\tilde{D}_{0} : = \frac{1}{m_{0}}(\tilde{\sigma}_{0}).
\]
We define the singular hermitian metric $\tilde{h}_{0}$ 
on $p^{*}K_{X}$ by 
\[
\tilde{h}_{0}:= \frac{1}{\mid \tilde{\sigma}_{0}\mid^{2/m_{0}}}.
\]
We shall replace $\alpha_{0}$  by 
\[
\tilde{\alpha}_{0} 
:= \inf \{\alpha > 0\mid 
\mbox{the generic point of}\,\, Z \subseteq \mbox{Spec}
({\cal O}_{X \times B}/{\cal I}(h_{0}^{\alpha}))\} .
\]
Then for every $0 < \delta << 1$, there exists a Zariski 
open subset $U$ of $B$ such that for every $b \in U$, 
$\tilde{h}_{0}\mid_{X\times\{ b\}}$ is well defined and  
\[
b \not{\subseteq}\mbox{Spec}({\cal O}_{X\times\{ b\}}/{\cal I}(\tilde{h}_{0}^{\alpha_{0}-\delta}\mid_{X\times\{ b\}})),
\]
where we have identified $b$ with distinct two points in $X$. 
And also by Lemma 2.6, we see that 
\[
b \subseteq \mbox{Spec}({\cal O}_{X\times\{ b\}}/{\cal I}(\tilde{h}_{0}^{\alpha_{0}}\mid_{X\times\{ b\}})),
\]
holds for every $b\in B$. 
Let $\tilde{X}_{1}$ be the minimal center of logcanonical singularities 
of $(X\times B,\alpha_{0}\tilde{D}_{0})$ at the generic point of $Z$. 
(although $\tilde{D}_{0}$ may not be effective this is meaningful 
by the construction of $\tilde{\sigma}_{0}$).   
We note that $\tilde{X}_{1}\cap q^{-1}(b)$ may not be 
irreducible even for a general $b\in B$. 
But if we take a suitable finite cover
\[
\phi_{0} : B_{0} \longrightarrow B,
\]
on the base change $X\times_{B}B_{0}$, $\tilde{X}_{1}$ 
defines a family of irreducible subvarieties
\[
f_{1} : \hat{X}_{1} \longrightarrow U_{0}
\]
of $X$ parametrized by a nonempty Zariski open subset
 $U_{0}$ of $\phi_{0}^{-1}(U)$.
Let $n_{1}$ be the relative dimension of $f_{1}$.
We set 
\[
\tilde{\mu}_{1} := K_{X}^{n_{1}}\cdot f_{1}^{-1}(b_{0})
\]
where $b_{0}$ is a general point on $U_{0}$.
Continueing this process 
we may construct a finite morphism 
\[
\phi_{r} : B_{r} \longrightarrow B
\]
and a nonempty Zariski open subset $U_{r}$ of $B_{r}$ 
which parametrizes a family of stratification 
\[
X \supset  X_{1} \supset X_{2} \supset \cdots \supset X_{r} 
\supset X_{r+1} =  x \,\,\mbox{or}\,\, x^{\prime}
\] 
constructed as before.
And we also obtain invariants $\{\tilde{\alpha}_{0},
\ldots ,\tilde{\alpha}_{r}\}$, $\{\tilde{\mu}_{0},\ldots ,\tilde{\mu}_{r}\}$,
$\{ n = \tilde{n}_{0}\ldots ,\tilde{n}_{r}\}$.
Hereafter we denote these invariants without $\,\,\tilde{} \,\,$ for simplicity.  By the same proof as Lemma 2.4, we have the following lemma. 
\begin{lemma}
\[
\alpha_{i}\leq \frac{n_{i}\sqrt[n_{i}]{2}}{\sqrt[n_{i}]{\mu_{i}}} + O(\varepsilon_{i-1})
\]
hold for $1\leq i\leq r$.
\end{lemma}
By Lemma 2.9 we obtain the following proposition. 
\begin{proposition}
For every 
\[
m \geq  \lceil\sum_{i=0}^{r}\alpha_{i}\rceil + 1
\]
$\mid mK_{X}\mid$ gives a birational rational map from $X$ into 
a projective space.
\end{proposition}

\subsection{Use of subadjunction theorem}
The following subadjunction theorem is crucial in our proof. 
\begin{theorem}(\cite{ka})
Let $X$ be a normal projective variety.
Let $D^{\circ}$ and $D$ be effective {\bf Q}-divisor on $X$ such that 
$D^{\circ} < D$, $(X,D^{\circ})$ is logterminal and 
$(X,D)$ is logcanonical. 
Let $W$ be a minimal center of logcanonical singularities for $(X,D)$. 
Let $H$ be an ample Cartier divisor on $X$ and $\epsilon$ a positive rational number.
Then there exists an effective {\bf Q}-divisor $D_{W}$ on $D$ such that 
\[
(K_{X}+D+\epsilon H)\mid_{W}\sim_{\mbox{\bf Q}}K_{W}+D_{W}
\]
and $(W,D_{W})$ is logterminal. 
In particular $W$ has only rational singularities.
\end{theorem}
\begin{remark}
As is stated in \cite[Remark 2.2]{ka3}, 
the assumption that $W$ is a minimal center can be replaced
that $W$ is a local minimal center, since the argument 
in \cite{ka} which uses the variation of Hodge structure 
does not change. 
But in this case we need to replace $K_{W}$ by the 
pushforward of the canonical divisor of the normalization 
of $W$, 
since \cite[p.494, Theorem 1.6]{ka2} 
works only locally in this case. 
\end{remark}
Roughly speaking, Theorem 2.2 implies that $K_{X}+D$ (almost) dominates $K_{W}$.  

Let $W_{j}$ be a nonsingular model of $X_{j}$. 
By Theorem 2.2, we see that 
\[
\mu (W_{j},K_{W_{j}})
\leq 
(1 +\sum_{i=0}^{j-1}\alpha_{i})^{n_{j}}\cdot \mu_{j}
\]
holds, where 
\[
\mu (W_{j},K_{W_{j}}) 
:= n_{j}!\cdot \overline{\lim}_{m\rightarrow\infty}
m^{-n_{j}}\dim H^{0}(W_{j},{\cal O}_{W_{j}}(mK_{W_{j}})).
\]
We note that if we take $x,x^{\prime}$ general, 
$W_{j}$ ought to be of general type. 
Otherwise $X$ is dominated by a family of varieties 
of nongeneral type. 
This contradicts the assumption that $X$ is of general type. 

For every smooth projective  variety $W$ of general type 
of dimension $k\leq 2$, 
\[
\mu (W,K_{W})\geq C(k)
\]
holds, where $C(1) = 2$ and $C(2) = 1$. 
Then by the above inequality  
\[
C(n_{j}) \leq 
(1 +\sum_{i=0}^{j-1}\alpha_{i})^{n_{j}}\cdot \mu_{j}
\]
holds. 
Since 
\[
\alpha_{i} \leq \frac{\sqrt[n_{i}]{2}n_{i}}{\sqrt[n_{i}]{\mu_{i}}}+ O(\varepsilon_{i-1})
\]
holds by Lemma 2.9, 
we see that  
\[
\frac{1}{\sqrt[n_{j}]{\mu_{j}}}\leq (1+\sum_{i=0}^{j-1}\frac{\sqrt[n_{i}]{2}n_{i}}{\sqrt[n_{i}]{\mu_{i}}})\cdot C(n_{j})^{-\frac{1}{n_{j}}}
\]
holds for every $j \geq 1$.
Using the above inequality inductively, we have
the following lemma.
\begin{lemma}
Suppose that $\mu_{0} \leq 1$ holds.
Then there exists a positive constant $C$ independent of $X$  
such that 
\[
\sum_{i=0}^{r}\frac{\sqrt[n_{i}]{2}n_{i}}{\sqrt[n_{i}]{\mu_{i}}}\leq \frac{C}{\sqrt[3]{\mu_{0}}}
\]
holds. 
In fact $C$ can be taken to be 
\[
C = 3\sqrt[3]{2}+ (3\sqrt[3]{2}+1)\sqrt{2} + 
\frac{1}{2} (3\sqrt[3]{2}+ (3\sqrt[3]{2}+1)\sqrt{2}+1). 
\]
\end{lemma} 
\subsection{Estimate of the degree}
\begin{lemma}
If $\Phi_{\mid mK_{X}\mid}\mid$ is birational rational map
onto its image, then
\[
\deg \Phi_{\mid mK_{X}\mid}(X)\leq \mu_{0}\cdot m^{3}
\]
holds.
\end{lemma}
{\bf Proof}.
Let $p : \tilde{X}\longrightarrow X$ be the resolution of 
the base locus of $\mid mK_{X}\mid$ and let 
\[
p^{*}\mid mK_{X}\mid = \mid P_{m}\mid + F_{m}
\]
be the decomposition into the free part $\mid P_{m}\mid$ 
and the fixed component $F_{m}$. 
We have
\[
\deg \Phi_{\mid mK_{X}\mid}(X) = P_{m}^{3}
\]
holds.
Then by the ring structure of $R(X,K_{X})$, we have an injection 
\[
H^{0}(\tilde{X},{\cal O}_{\tilde{X}}(\nu P_{m}))\rightarrow 
H^{0}(X,{\cal O}_{X}(m\nu K_{X}))
\]
for every $\nu\geq 1$.
We note that since ${\cal O}_{\tilde{X}}(\nu P_{m})$ is globally generated
on $\tilde{X}$, for every $\nu \geq 1$ we have the injection 
\[
{\cal O}_{\tilde{X}}(\nu P_{m})\rightarrow p^{*}({\cal O}_{X}(m\nu K_{X}).
\]
Hence there exists a natural morphism 
\[
H^{0}(\tilde{X},{\cal O}_{\tilde{X}}(\nu P_{m}))
\rightarrow 
H^{0}(X,{\cal O}_{X}(m\nu K_{X}))
\]
for every $\nu\geq 1$. 
This morphism is clearly injective. 
This implies that 
\[
\mu_{0} \geq  m^{-3}\mu (\tilde{X}_{i},P_{m})
\]
holds. 
Since $P_{m}$ is nef and big on $X$ we see that 
\[
\mu (\tilde{X},P_{m}) = P_{m}^{3}
\]
holds.
Hence
\[
\mu_{0}\geq m^{-3}P_{m}^{3}
\]
holds.  This implies that
\[
\deg \Phi_{\mid mK_{X}\mid}(X)\leq \mu_{0}\cdot m^{3}
\]
holds.
{\bf Q.E.D.}
\vspace{5mm} \\
By Lemma 2.9.2.10,2.11 we see that 
if $\mu_{0} \leq 1$ holds, 
for 
\[
m := 1+\lceil \sum_{i=0}^{r}\frac{\sqrt[n_{i}]{2}n_{i}}{\sqrt[n_{i}]{\mu_{i}}}\rceil 
\]
$\mid mK_{X}\mid$ gives a birational embedding of $X$ and 
\[
\deg \Phi_{\mid mK_{X}\mid}(X) 
\leq C^{3}
\]
holds. 
Since $C$ is a positive constant independent of $X$, 
we have that there exists a positive constant $C(3)$ 
independent of $X$ such that 
\[
\mu_{0} = K_{X}^{3} \geq C(3)
\]
holds. 

More precisely we arugue as follows. 
Let ${\cal H}$ be an irreducible component of the Hilbert scheme of a projective space. 
Let ${\cal H}_{0}$ be the Zariski open subset of ${\cal H}$ which parametrizes 
irreducible subvarieties. 
Then there exists a finite stratification of ${\cal H}_{0}$ by Zariski locally closed subsets such that on each strata there exists a simultaeneous resolution 
of the universal family on the strata. 
We note that the volume of the canonical bundle of the resolution is constant on each strata by \cite{tu6,nak}.
Hence there exists a positive constant $C(n)$ such that 
\[
\mu (X,K_{X})\geq C(3)
\]
holds for every projective 3-fold $X$ of general type.

Then by Lemma 2.9 and 2.10, we see that there exists 
a positive integer $\nu_{3}$  such that 
for every projective $3$-fold $X$ of general type, 
$\mid mK_{X}\mid$ gives a birational embedding into a 
projective space for every $m\geq \nu_{3}$, 
This completes the proof of Theorem 2.1.

\section{Bounding the index of $X$}
Let $X$ be a canonical 3-fold of general type such that 
$\mu (X,K_{X})= K_{X}^{3}$ is minimal among all 
canonical 3-fold of general type. 
\subsection{3-dimensional terminal singularities}
Let $r,a_{1},a_{2},a_{3}$ be coprime positive integers. 
Let $\xi$ be a primitive $r$-th root of unity acting 
on $\mbox{\bf C}^{3}$ via 
\[
\xi (x,y,z) = (\xi^{a_{1}}x,\xi^{a_{1}}y,\xi^{a_{3}}z).
\]
A singularity $Q\in X$ is of type $\frac{1}{r}(a_{1},a_{2},a_{3})$,
if $(X,Q)$ is isomorphic to 
$(\mbox{\bf C}^{3},O)/\langle \xi\rangle$.
It is known that $\frac{1}{r}(a_{1},a_{2},a_{3})$
is terminal, if and only if 
$(a_{1},a_{2},a_{3}) = (1,-1,b)$. 

In general a 3-dimensional terminal singularity is not 
necessarily a quotient singularity. 
But we can describe a 3-dimensional terminal singularity as 
follows.

Let $G_{r}$ be the cyclic group of order $r$.
Let $\xi$ be a generator of $G_{r}$. 
Let 
\[
\tilde{X} = \{ (x,y,z,t)\in \mbox{\bf C}^{4}\mid 
xy- g(z^{r},t)\} ,
\]
where $g$ is a  polynomial. 
Let us define the action of $G_{r}$ on $\tilde{X}$ 
by 
\[
\xi (x,y,z,t) := (\xi x,\xi^{-1}y,\xi^{b}z,t),
\]
where $b$ is a positive integer coprime to $r$. 
We set $X := \tilde{X}/G_{r}$.
We may consider $X$ as a family over {\bf C} by 
the projection $(x,y,z,t)\mapsto t$. 
By \cite{m}, all the 3-dimensional terminal singularities
of index $r\geq 5$ are  obtained as a member of 
such a family $X \longrightarrow {\bf C}$. 
We note that since $X$ is a Cartier divisor 
in $\frac{(1,-1,b)}{r}\times \mbox{\bf C}$,
all the 3-dimensional terminal singularities 
of index $r\geq 5$ can be deformed to 
a sum of quotient singularities of type 
$\frac{(1,-1,b)}{r}$.
The 3-dimensional terminal singularities of index $\leq 4$ 
are described similary. 

\begin{lemma}
Let $(M,Q)$ be a germ of a 3-dimensional terminal singularity 
of index $r \geq 5$. 
Let $\pi : (N,P)\longrightarrow (M,Q)$ be the canonical 
cover, i.e. $\pi$ is a cycric $r$-covering associated 
with the isomorphim ${\cal O}_{M}(rK_{M})\simeq {\cal O}_{M}$. 
Suppose that $(M,Q)$ can be deformed 
to a sum of $\frac{(1,-1,b)}{r}$ as above. 
Then for a positive integer $m$, 
\[
\dim H^{0}(N,
({\cal O}_{N}(mK_{N})/\pi^{*}{\cal O}_{M}(mK_{M}))^{\otimes \ell})
\geq  \frac{1}{6}\overline{(r-m)}\cdot\overline{(r-bm)}\cdot\overline{bm}\cdot \ell^{3}
+ O(\ell^{2})
\]
holds, where for an integer $a$, $\bar{a}$ denotes 
the smallest nonnegative residue modulo $r$.  
\end{lemma}
{\bf Proof}. 
If $(M,Q)$ is $\frac{(1,-1,b)}{r}$, 
the lefthandside is equal to the righthandside by 
an easy direct calculation (just by counting 
invariant monomials under the action of $G_{r}$). 
Otherwise since $(M,Q)$ is a hypersurface in 
$\frac{(1,-1,b)}{r}\times \mbox{\bf C}$ as above, 
the assertion is again clear. 
{\bf Q.E.D.}

\subsection{Boundedness of the degree}

We set $\mu_{0}:= K_{X}^{3}$. 
Suppose that that $\mu_{0} \leq 1$ holds. 
Let us set 
\[
m_{0}:= \lceil 18\frac{1}{\sqrt[3]{\mu_{0}}} \rceil
\geq 
\lceil (3\sqrt[3]{2}+1)+ (3\sqrt[3]{2}+1)\sqrt{2} + 
\frac{1}{2} (1+3\sqrt[3]{2}+ (3\sqrt[3]{2}+1)\sqrt{2}+1))
\frac{1}{\sqrt[3]{\mu_{0}}} \rceil .
\]
Then by Lemma 2.11,
\[
\deg \Phi_{\mid m_{0}K_{X}\mid}
\leq \mu_{0}\cdot m_{0}^{3} \leq 18^{3}
\]
hold.

By Lemma 3.1 and the proof of Lemma 2.11, 
we have that  
\[
\deg_{\mid m_{0}K_{X}\mid}(X)\leq \mu_{0}m_{0}^{3} -
\sum_{Q}\frac{\overline{(r-m_{0})}\cdot \overline{bm_{0}}\cdot\overline{(r-bm_{0})}}{6r}
\]
holds, where $Q$ runs the basket of singularities on $X$. 

We set $r_{0} = \max_{Q} r$. 
Suppose that 
$r \geq 2m_{0}$ holds. 
Then we have that 
\[
\frac{\overline{(r_{0}-m_{0})}\cdot\overline{bm_{0}}\cdot\overline{(r_{0}-bm_{0})}}{6r_{0}} 
\geq \frac{1}{12}\overline{bm_{0}}\cdot\overline{(r_{0}-bm_{0})}
\geq \frac{r_{0}-1}{12}
\]
hold. 

Suppose that $r_{0}\leq 2m_{0}$ holds. 
In this case, there exists a positive integer 
$m_{0}\geq m^{\prime}_{0}\leq 2m_{0}$ such that 
\[
\frac{\overline{(r_{0}-m_{0}^{\prime})}\cdot\overline{bm_{0}^{\prime}}\cdot\overline{(r_{0}-bm^{\prime}_{0})}}{6r_{0}}
\geq \frac{\overline{(r_{0}-m_{0}^{\prime})}}{6r_{0}}\cdot
\frac{r_{0}^{2}}{4} 
\geq \frac{r_{0}}{24}
\]
hold. 
Hence combining these, we have that 
\[
r_{0}\leq 2^{6}\cdot 3\cdot 18^{3} = 2^{9}\cdot 3^{7}
\]
hold. 

\subsection{A bound for $\nu_{3}$}
In the above situation, the global index of $X$, i.e.,
the least positive integer $R$ such that $RK_{X}$ is Cartier
is less than  
\[
r_{0}! \leq (2^{9}\cdot 3^{7})!.
\]
This implies that 
\[
\mu_{0} \geq \frac{1}{((2^{9}\cdot 3^{7})!)^{3}}
\]
Hence by the formulas at the beggining of Section 3.2, we see that 
\[
m_{0}\leq  18(2^{9}\cdot 3^{7})!
\]
holds. 
This completes the proof of Theorem 1.1.

Author's address\\
Hajime Tsuji\\
Department of Mathematics\\
Tokyo Institute of Technology\\
2-12-1 Ohokayama, Megro 152-8551\\
Japan \\
e-mail address: tsuji@math.titech.ac.jp

\end{document}